\documentclass[10pt]{article}
\font\smallit=cmti10

\usepackage{float}
\usepackage{amssymb,amsmath,amsthm,enumitem} 
\usepackage{subcaption,booktabs}
\usepackage{hyperref}
\usepackage{tikz}
\usepackage{tabularray} 
\UseTblrLibrary{booktabs}
\usetikzlibrary{positioning,shapes.misc,shapes.geometric, graphs, calc}
\makeatletter
\renewcommand\section{\@startsection {section}{1}{\z@}
{-30pt \@plus -1ex \@minus -.2ex}
{2.3ex \@plus.2ex}
{\normalfont\normalsize\bfseries\boldmath}}

\renewcommand\subsection{\@startsection{subsection}{2}{\z@}
{-3.25ex\@plus -1ex \@minus -.2ex}
{1.5ex \@plus .2ex}
{\normalfont\normalsize\bfseries\boldmath}}

\renewcommand{\@seccntformat}[1]{\csname the#1\endcsname. }

\makeatother

\newtheorem{theorem}{Theorem}

\newtheorem{cnj}{Conjecture}
\newtheorem{proposition}{Proposition}

\newtheorem{example}{Example}

\theoremstyle{definition}

\newtheorem{remark}{Remark}

\def\gO{\Omega}
\def\go{\omega}

\def\cG{\mathcal{G}}

\begin{document}
\begin{center}
\uppercase{A counterexample to conjecture ``Catch 22"}
%with 3 players, and 5 outcomes: 2 terminal and 3 cyclic}

\vskip 10pt
{\bf Bogdan Butyrin}\\
{\smallit National Research University 
Higher School of Economics (HSE), Moscow, Russia}\\
{\tt butirin.bogdan@yandex.ru}\\
\vskip 5pt
{\bf Vladimir Gurvich}\\
{\smallit National Research University Higher School of Economics (HSE), Moscow, Russia}\\
{\tt vgurvich@hse.ru}, {\tt vladimir.gurvich@gmail.com}\\
\vskip 5pt
{\bf Anton Lutsenko}\\
{\smallit National Research University Higher School of Economics (HSE), Moscow, Russia}\\
{\tt ailutsenko@edu.hse.ru}\\
\vskip 5pt
{\bf Mariya Naumova}\\
{\smallit Rutgers Business School, Rutgers University, Piscataway, NJ, United States}\\
{\tt mnaumova@business.rutgers.edu}\\
\vskip 5pt
{\bf Maxim Peskin}\\
{\smallit National Research University Higher School of Economics (HSE), 
Moscow, Russia}\\
{\tt mppeskin\_1@edu.hse.ru}\\
\end{center}
\vskip 10pt
%\centerline{\smallit Received: , Revised: , Accepted: , Published: }
% We will fill in the dates
\vskip 30pt

{\bf Abstract:}
%\centerline{\bf Abstract} \noindent
We construct a finite deterministic graphical (DG) game 
without Nash equilibria  
in pure stationary strategies. 
This game has 3 players $I=\{1,2,3\}$  and 5 outcomes:   
2 terminal $a_1$ and $a_2$ and 3 cyclic. 
Furthermore, for 2 players a terminal outcome is the best:  
$a_1$ for player 3 and $a_2$ for player 1. 
Hence, the rank vector $r$ is at most $(1,2,1)$.  
Here $r_i$ is the number of terminal outcomes 
that are worse than some cyclic outcome 
for the player $i \in I$.  
This is a counterexample to 
conjecture ``Catch 22" from the paper 
``On Nash-solvability of finite $n$-person DG games, Catch 22" (2021)  
\url{https://arxiv.org/abs/2111.06278}, 
according to which,   
at least 2 entries of $r$ are at least 2 for any NE-free game. 
However, Catch 22 remains still open for 
the games with a unique cyclic outcome, 
not to mention a weaker 
(and more important) conjecture 
claiming that an $n$-person finite DG 
game has a Nash equilibrium 
(in pure stationary strategies) 
when $r = (0^n)$, that is, all $n$ entries of $r$ are 0;  
in other words, when the following condition holds: 
\newline 
($C_0$)  
any terminal outcome is better than every cyclic 
one for each player. 
\newline 
A game is play-once 
if each player controls a unique position. 
It is known that any play-once game 
satisfying ($C_0$) has a Nash equilibrium.
We give a new and very short proof of this statement.
Yet, not only conjunction but already disjunction 
of the above two conditions may be sufficient for 
Nash-solvability. This is still open.
\newline
{\bf Keywords:}
$n$-person deterministic graphical %(multi-stage) 
games,
Nash equilibrium, %Nash-solvability,
pure stationary strategy, digraph, directed cycle, 
strongly connected component.
%\end{abstract}

\section{Introduction; main conjectures and examples}
\label{s0}
Our main example is given 
in graphical and normal forms in Figures \ref{f1} and \ref{f2}, 
respectively. 
Here we briefly describe its main properties 
postponing basic definitions, 
which will be given later in Introduction.
Figure \ref{f1} represents a finite 
deterministic graphical game of 3 players, 
$I = \{1,2,3\}$,  controlling positions 
$V_1 = \{s_1,u_1,v_1\}; \; V_2 = \{u_2,v_2\}$ and 
$V_3 = \{u_3,v_3\}$, respectively. 

We restrict ourselves to finite games 
with perfect information and without moves of chance, 
and the players to their pure stationary strategies. 
In our example, player 1 has $3 \times 2 \times 2 = 12$ strategies, 
while players 2 and 3 have $2 \times 2 = 4$ strategies, each. 
This game in the normal form is defined by 
the table $X = X_1 \times X_2 \times X_3$ 
of size $12 \times 4 \times 4$,    
the game form mappings $g : X \rightarrow \gO$, 
and a triple of preference relations   
$\succ = (\succ_1, \succ_2, \succ_3)$  over $\gO$,   
where $\gO = \{a_1, a_2, c_1, c_2, c_3\}$ 
is the set of outcomes. 
It consists of  5  strongly connected components 
(SCCs) of the directed graph (digraph) $G = (V,E)$ of the game. 
Two of them are the terminals: $a_1$ and $a_2$, 
while 3 are cycles: $c_1, c_2,$ and $c_3$; see Figure \ref{f1}. 
Note that there is one more SCC in $G$, 
it consists of the initial position $s_1$. 
Such an SCC, which is not terminal and 
does not contain a directed cycle, 
is called {\em transient}; 
no outcome is assigned to it.

A digraph $G = (V,E)$ is partitioned into SCCs. 
This partition is unique
and can be constructed 
in time linear in its size of $G$, 
that is, $(|V | + |E|)$. 
This partition has many important properties and applications; 
see, for example,\cite{Tar72} and \cite{Sha81} for  details.
By contracting each SCC of $G$ to a single vertex 
(and deleting every edge that belongs to an SCC)  
one obtains an acyclic digraph $G^*$. 
One more application, to the DG games, was suggested in \cite{Gur18}.
The main result of this paper is the Nash-solvability 
(in pure stationary strategies) 
of the (finite) $2$-person ($n=2$) DG games.  
It follows from the general result 
on Nash-solvability of tight game forms 
\cite{Gur75,Gur89}; see also \cite{GN21,GN21A,GN22}.
For the zero-sum case, the saddle-point-solvability 
was proven even earlier \cite{EF70}, see also \cite{Gur73}. 

However, this result cannot be extended to $n>2$. 
The first example of an NE-free $4$-person DG game  
was constructed in \cite{GO14}; see also \cite{Gur15}. 
Then a much simpler $3$-person NE-free DG game 
was given in \cite[Figure 3]{BGMOV18}. 
It has 4 outcomes: 3 terminal and a unique cycle. 
Furthermore, its digraph $G$ has two transient SCCs, 
which are not associated with outcomes. 
Moreover, this game is ``almost play-once": 
players 1 and 3 control a unique position each,  
only player 2 controls two positions.  
%According to \cite{BG03} a game is called 
%{\em play-once} if each player controls a unique position. 

The remaining important part of a DG game 
is the payoff $u : I \times \gO \rightarrow \mathbb{R}$, 
where $u(i, \omega)$  is the profit of player $i \in I$ 
in case outcome $\go \in \gO$ is realized in the game. 
Studying NE, it is both correct and convenient  
to replace payoffs by preferences. 
For each player $i \in I = \{1, \dots, n\}$, 
we introduce a linear order $\succ_i$  over $\Omega$. 
Notation $\go' \succ_i \go''$  means that 
$i$  prefers $\go'$ to $\go''$.
In particular, we assume that there are no ties. 
It can be done without any loss of generality (wlog). 
Indeed, making ties can create an NE, but cannot destroy it; 
see the definition of Subsection \ref{ss-NE}. 
% later in Introduction. 
In particular, an NE-free  DG game with payoff $u$ 
remains NE-free after any sufficiently 
small perturbation of $u$. 
Thus, one can replace a payoff $u$ 
by a close tie-free payoff $u'$ and 
then, in its turn, replace $u'$
by the corresponding $n$-dimensional preference profile 
$\succ = \{\succ_i \mid i \in I\}$, 
keeping the transformed game NE-free all way.  

Note also that merging arbitrarily the outcomes 
is an operation that respects NE-existence; 
see for example, \cite{Gur75,Gur89}. 
In this paper, merging all cyclic outcomes into a unique one, 
will play an important role. 

In \cite{BG03}, among other results, it was shown that 
a play-once DG game has an NE 
if the following additional condition hold: 

\medskip 

CND$(C_0)$: \; 
every terminal outcome is better 
than any cyclic one for each player. 

\medskip 

A very simple proof of this statement 
will be given in Subsection \ref{ss-01}. 

Clearly, CND$(C_0)$ holds if and only if 
the rank $n$-vector $r$ is $(0^n)$, 
that is, all $n$ entries of $r$  equal 0.

However, it seems that conjunction of 
the above two conditions:
the play-once one and $(C_0)$ is too strong. 
Can it be replaced by disjunction? 
No counterexample is known, yet. 

In \cite[Section 6]{BGMOV18}, it was shown that 
if an NE-free DG game satisfying  CND$(C_0)$ exists, 
then there is also a DG game that has  
a unique NE-outcome, which is cyclic. 
This conclusion holds even if  
the original game has terminal outcomes 
and only one cyclic one.
This looks like a paradox in view of $(C_0)$. 

Although CND$(C_0)$ looks ``relaxable", 
all attempts to construct an  
NE-free example satisfying $(C_0)$ fail yet. 
The following conjecture was suggested in \cite{Gur21}: 

\medskip 

CNJ$(C_0)$: \; 
A DG game has an NE whenever CND$(C_0)$ holds.

\medskip 

Let us merge all cyclic outcomes in one.
Then, CND$(C_0)$ (resp.  CNJ$(C_0)$) 
is replaced by a stronger condition CND$(C'_0)$ 
(resp., by a weaker conjecture CNJ$(C'_0))$. 

Since attempts to disprove 
the stronger CNJ$(C'_0)$ fail too, 
we try to strengthen it further. 
Given the preference $\succ_i$ 
of a player $i \in I$,  let  
$r_i$  be the number of terminal outcomes for which  
there exists a better, for player  $i$, cyclic outcome.  
Thus every preference profile 
$\succ = \{\succ_i \mid i \in I\}$  
uniquely defines the rank $n$-vector 
$r = (r_i \mid i \in I)$. 
For example, for 
the DG game from \cite{BGMOV18} we have $r = (2,0,2)$. 
For many other NE-free DG games their rank vector 
has at least 2 entries equal to at least 2, each. 

In \cite{Gur21}, 
the corresponding condition and conjecture 
were denoted by 
CND$(C_{22})$ and CNJ$(C_{22})$, respectively. 
Furthermore, after merging all cyclic outcomes in one, 
we obtain CND$(C'_{22})$ and CNJ$(C'_{22})$. 
CNJ$(C'_{22})$, which is stronger than CNJ$(C_{22})$, 
holds for NE-free DG games from \cite{GO14,BGMOV18} 
and in many other cases, as our computer analysis shows. 
Partial results were obtained in \cite{GN22}.

\medskip 

Our main example disproves CNJ$(C_{22})$ but not CNJ$(C'_{22})$. 
Indeed, consider the preference profile, given in Figure \ref{f1}. 
There are 3 cyclic and only 2 terminal outcomes. 
Furthermore, one of them is the best outcome, for 2 players. 
In other words,  $r = (1,2,1)$  % [MAYBE A PERMUTATION ?]. 
However, an NE appears after merging $c_1,c_2,$ and $c_3$.  

\medskip 

To summarize, CNJ$(C_{22})$ is disproved, 
but CNJ$(C'_{22})$ still stays, along with 
CNJ$(C_0)$ and CNJ$(C'_0)$. 

\medskip 

A digraph is called symmetric if 
$(u,v)$ is its directed edge whenever $(v,u)$ is 
(unless $u$ is a terminal or $v$ is an initial position).

In \cite[Theorem 2]{BFGV23}, it was shown that 
every $n$-person DG game on a symmetric digraph is Nash-solvable; 
see also Theorem 3.

\bigskip 

We assume, yet, that all cycles of any fixed SCC 
form the same outcome. 
Another option is to wave this assumption 
and consider  each cycle as a separate outcome. 
Then, Nash-solvability becomes a rare phenomenon 
already for the 2-person zero-sum games. 
Consider, for example, the digraph with 
two vertices $s,v$ controlled by players 
1 and 2, respectively, 
and 2 pairs of directed edges: 
two from $s$ to $v$ and two from $v$ to $s$;
its normal form is the $2 \times 2$  
matrix game whose entries are 
four distinct outcomes; see \cite[Figure 1 (1)]{BGMW11}.
Thus, no saddle point exists  
if player 1 prefers the main diagonal to the side one, 
while player 2 prefers vice versa.

A criterion of Nash-solvability 
for such two-person DG games on symmetric digraphs 
is given by \cite[Theorem 1, Corollaries 1 and 2]{BGMW11}. 

\section{Main concepts and definitions}
\label{s1}
\subsection{Strongly connected components  
of directed graphs}
\label{ss10}

Let $G = (V,E)$ be a directed graph (digraph). 
It may have parallel edges, 
directed oppositely or similarly, and loops, 
at most one at each vertex. 
Digraph $G$ is called {\em strongly connected} 
if for any $v, v' \in V$ there is a directed path from 
$v$ to $v'$ (and, hence, from $v'$ to $v$, as well). 
By definition, the union of two strongly connected digraphs 
is strongly connected if and only if they have a common vertex. 
A vertex-inclusion-maximal strongly connected induced subgraph of G 
is called its {\em strongly connected component} (SCC). 
Obviously, any digraph $G = (V, E)$ 
admits a unique decomposition into SCCs:
$G^j = G[V^j] = (V^j, E^j)$ for $j \in J$, 
where $J$ is a set of indices and 
$V = \cup_{j \in J} V^j$  is a partition of $V$, 
that is, $V^j \cap V^{j'} = \emptyset$ 
whenever $j \neq j'$  for a pair $j, j' \in J$. 
This partition can be determined in time linear 
in the size of $G$ (that is, in $(|V |+|E|))$ 
and has numerous applications; 
see \cite{Tar72} and also \cite{Sha81} for more details. 
One more application,  
to positional games, was suggested in \cite{Gur18}.

Let $V_T \subseteq V$ 
denote the set of all the terminals of $G$, 
that is, $v \in V_T$ if and only if there is no edge from $v$.
% , except, perhaps, a loop $e_v = (v, v)$. 
Any such vertex forms a maximal SCC, that is, $\{v\} = V^j$ 
for some $j \in J$. 
Denote by $J^T$ the corresponding subset of $J$. 
For convenience, let us add a loop $e_v$ to each $v \in V_T$ 
that does not have one yet. All loops in $V_T$ , 
old and new, will be called terminal loops and all 
$v \in V_T$ will be still called terminals.

The concepts of a directed path and a directed cycle 
(dicycle) are defined standardly \cite{Sha81,Tar72}. 
An edge is a directed path of length 1; 
a loop is a dicycle of length 1, 
a dicycle of length 2 is a pair of parallel oppositely directed edges, 
a dicycle of length $|V|$ 
is called {\em Hamiltonian}. 
A digraph without dicycles is called {\em acyclic}.

Let $J_0$ denote the subset of $J$ such that 
$V^j$ is a single vertex with no loop it it. 
Obviously, $G^j$ contains a dicycle unless $j \in J_0$. 
By our convention, $J_T \cap J_0 = \emptyset$.

For each $j \in J$, let us contract $G^j$ 
into a single vertex $v^j$. 
Then, all the edges of $E^j$ 
(in particular, loops) disappear and
we obtain an {\em acyclic} digraph  $G^* = (V^*,E^*)$. 
It can be constructed in time linear in $|V|+|E|$. 

\subsection{Modeling graphical games by digraphs} 
\label{ss00}
Given an arbitrary digraph $G = (V, E)$, 
we interpret it as a {\em positional structure}  
in which $V$ and $E$ are positions and moves, respectively. 
Furthermore, let us fix a partition 
$D : V = V_1 \cup \dots \cup V_n \cup V_T$ 
in which $V_i$ are the sets of positions controlled by the
players $i \in I = \{1, \dots, n\}$, 
while $V_T$ is the set of all terminals.

The considered model is restricted to the games 
with perfect information and without moves of chance. 
Also, we restrict
all players $i \in I$ to their 
pure stationary strategies $X_i$.  
Such a strategy $x_i \in X_i$ assigns a move
$(v, v') \in E$ to every $v \in V_i$. 
Given an initial position $v_0 \in V \setminus V_T$, 
each {\em strategy profile} $x = (x_1, \dots, x_n)$ 
uniquely defines a walk
in $G$ that begins in $v_0$ and 
in each position follows the move chosen 
by the corresponding strategy. 
We denote this walk by $p(x)$ and 
call it a {\em play} generated by $x$.  
Since all strategies of $x$ are stationary, 
play $p(x)$  forms a “lasso”, that is, 
it consists of an initial directed path, 
which can be empty, and 
a dicycle repeated infinitely. 

For every $j \in J \setminus J_0$, 
let us identify all lassos whose dicycles belong to 
$G_j$ and treat them as a single outcome $c_j$ of the game. 
Then, $\gO = \{c_j \mid j \in J \setminus J_0\}$ 
is the set of outcomes. 

An outcome is called {\em terminal} 
if the corresponding SCC $G^j$ 
generates a terminal vertex of $G^*$. 
Note that all terminals of  $G$  
are the terminals $G^*$, but there may exist others.   
Any SCC such that each move from it 
does not leave it, also forms a terminal 
of the considered game. 
Furthermore, components $G^j$  for $j \in J_0$ 
are called {\em transient}. 
Each one consists of a single vertex and 
does not form an outcome. 
Each remaining SCC, not terminal and not transient, 
will be called {\em inner}. 
Clearly, each one has a cycle, which may be a loop, 
and it forms an outcome, 
which will be called {\em inner} too. 

Thus, given a digraph 
$G = (V,E)$  and partition $D$,  
we obtain mapping  $g : X \rightarrow \gO$, 
where $X = X_1 \times \dots \times X_n$  
is the set of strategy profiles, also called 
{\em situations}. 
In general, such a mapping $g$ is called a {\em game form}. 

Pair $(G,D)$ defines a {\em DG game structure} 
and $g = g(G,D)$  its {\em normal form}. 

Given an $n$-person game form 
$g : X \rightarrow \gO$, introduce a payoff 
$u : I \times \gO \rightarrow \mathbb{R}$. 
The value $u(i,\go)$  is interpreted 
as the profit of the player $i \in I = \{1, \dots, n\}$ 
in case the outcome $\go \in \gO$  is realized. 
Pair $(g,u)$ is called the {\em game in normal form}. 

The DG games were introduced by Washburn \cite{Was90} 
in 1990 and since then are frequent in the literature; 
see, for example,  
\cite{AGH10,AHMS12,BEGM12,BG03,%[Section 12]
BGMW11,Con92,DS02,GO14}.
Only 2-person zero-sum DG games were considered in 
\cite{Was90} %,Con92,AHMS12}, 
but the extension to the $n$-person case is obvious. 

\subsection{Nash equilibrium (NE)} 
\label{ss-NE}
Given the game $(g,u)$ defined 
by a game form $g : X \rightarrow \gO$  
and payoff $u : I \times \gO \rightarrow \mathbb{R}$, 
a situation 
$x = (x_1, \dots, x_n) \in X_1 \times \dots X_n = X$ 
is called a {\em NE} if for each player $i \in I$ 
inequality $u(i, g(x)) \geq u(i, g(x'))$  
holds for every  situation  $x' \in X$  
that may differ from $x$ only in entry $i$, that is, 
equality  $x_j = x'_j$  holds for all $j \neq i$.
In other words, no player $i \in I$ can  
improve the result  $u(i, g(x))$  
replacing the strategy  $x_i$  by $x'_i$  
provided all other players keep their strategies unchanged. 
This concept was introduced 
by John Nash \cite{Nas50,Nas51}. 
He proved that any matrix game has a NE in mixed strategies. 
Here we study Nash-solvability in pure strategies; 
see, for example,  
%\cite{ARV09,AGH10,AHMS12,BEGM12,BFGV23,BG03,BG09,BGMOV18,BGMW11,EF70,Gur73,Gur75,Gur88,Gur89,Gur15,Gur17,Gur18,Gur21,Gur21a,GK18,GN21,GN21A,GN21B,GN22,GO14,HL97,HR04,KMMM90,Mil96,Mil96_2,MS96,Ros73,Was90}.
[1-9,12-33,36,40].

%As we already mentioned, 
%tie-breaking may destroy an NE but cannot create it. 
%More precisely, replace a payoff $u$  by $u'$  such that 
%$u'(\go) \geq u'(\go')$  whenever $u(\go) \geq u(\go')$ 
%for any $\go, \go' \in \gO$. 
%(Yet, equality $u'(\go) = u'(\go')$ may replace 
%inequality $u(\go) > u(\go')$.)  
%By definition, $x$ is an NE in game $(g,u')$ 
%if it was an NE in game $(g,u)$ (but not vice versa). 

A game form $g$ is called {\em Nash-solvable} 
if for an arbitrary payoff  $u$, 
the corresponding game $(g,u)$ has an NE.  

Since tie-breaking can only destroy NE, but 
cannot create one, studying Nash-solvability, 
we can assume wlog that payoff  $u$  has no ties
and then replace $u$ 
by a preference $\succ = (\succ_1, \dots, \succ_n)$ 
assuming that $\go \succ_i \go'$  
if and only if  $u(i, \go) > u(i, \go')$ 
for any $\go, \go' \in \gO$. 

Recall also that merging outcomes 
of a game form respects its Nash-solvability. 

\section{Partial results}
\label{s2}
\subsection{Nash-solvability of $n$-person DG games
for $n=2$ and $n>2$} 
\label{ss20}

\begin{theorem} (\cite{Gur18}) 
\label{t-n=2}
Every 2-person DG game has an NE.
\end{theorem}

This result is derived from a general 
criterion of Nash-solvability of 2-person game forms 
\cite{Gur75,Gur89}. 
The zero-sum case was considered even earlier 
by Edmonds and Fulkerson \cite{EF70}; see also \cite{Gur73}.
In \cite{Gur89} it was shown that 
the zero-sum- and Nash-solvability are equivalent properties 
of a 2-person game form $g$. 
Furthermore, both are equivalent with the tightness of $g$.
Tightness of the DG game forms was shown 
in \cite{BG03} for the special case, 
when all inner cyclic outcomes are merged, 
and verified in general in \cite{Gur18}. 
However, tightness is not related 
(neither necessary, nor sufficient) 
to Nash-sovability of $n$-person game forms 
with $n>2$ \cite{Gur89}.
Moreover, an $n$-person DG game is tight for all $n$, 
yet, may be not Nash-solvable when $n>2$; 
see \cite{GO14} for $n=4$ and \cite{BGMOV18} for $n=3$. 

NE-free DG games are sparse. 
We can add some extra conditions and enforce Nash-solvability. 
For example, it holds for all $n$ in case of
symmetric digraphs \cite{BFGV23}. 

Extra conditions 
CND$(C_0)$, CND$(C'_0)$, CND$(C_{22})$, CND$(C'_{22})$ 
and the corresponding conjectures on Nash-solvability, 
CNJ$(C_0)$, CNJ$(C'_0)$, CNJ$(C_{22})$, CNJ$(C'_{22})$, 
were considered in \cite{Gur21}. 
Here we disprove  CNJ$(C_{22})$, 
while 3 others remain open. 

\subsection{Play-once $n$-person DG games 
and condition CND$(C_0)$}
\label{ss-01}

The following two conditions are sufficient for 
the existence of a NE.

\begin{proposition}
An $n$-person DG game has an NE if  
(i) it is play-once and (ii) CND$(C_0)$ holds. 
\end{proposition}

This statement follows from the results of \cite{BG03}. 
Here we give a simple direct proof. 
Let us note, however, that not only 
conjunction but, perhaps, already 
disjunction of (i) and (ii) 
is sufficient for Nash-solvability. 
At least, we have no counterexample. 

So conjecture CNJ($C_{22}$) is disproved,  
but to conjectures 
CNJ($C_0$), CNJ($C'_0$) and CNJ($C'_{22}$) 
we can add one more. 

\begin{cnj}
A play-once $n$-person DG game is Nash-solvable.
\end{cnj}

Although a NE exists whenever 
conditions 
(i) and (ii) holds, yet, improvement cycles 
can exist too; see \cite[Figure 5 and Section 2]{AGH10} 
for definitions and more details. 

\medskip 

{\em Proof of Proposition.}
Note first that each 
(finite) digraph $G$  has a terminal SCC, 
since its SCC digraph $G^*$ is finite and acyclic. 

Assume for contradiction that 
a DG game $\cG = (G,D,\succ)$ is NE-free 
but satisfies (i) and (ii). 
Moreover, assume that it is a minimal such game. 
Wlog, assume that 
every vertex (and arc) of $G$ is reachable 
from the initial position $s$. 
Indeed, one can simply delete from 
$G$ the vertices not reachable from $s$, 
thus, reducing $G$.

{\bf Deleting bad terminal moves.} 
Suppose there are two terminal moves 
$v,t$  and  $v,t'$. 
Wlog we assume that 
$v$ is controlled by a player  $i \in I$ 
and that  $t \succ_i t'$. 
Then, one can delete arc $(v, t')$ 
and, obviously, no NE appears. 
Thus, we can assume that in every position  $v$  
there exists at most one terminal move. 
Furthermore, clearly, 
there is at least one non-terminal move in $v$. 
More generally, $\cG$ has no forced move. 
Indeed, in contradiction with  minimality of $G$, 
one can contract 
such a move keeping the normal normal form of the game 
and, in particular, the set of its Nash equilibria. 

\medskip 

Choose a shortest terminal play  $p$  
and fix the following situation: 
on $p$  follow  $p$, 
out of $p$ choose a non-terminal move. 
By (i) and (ii), this is a NE situation. 
\qed 

\begin{remark}
Recall that each finite DG game 
has a terminal and a terminal play.
Indeed, the game may have no terminal positions 
but it always has a terminal SCC, 
since its SCC graph is finite and acyclic.  
In the original graph, one can contract 
a terminal SCC replacing it by a single vertex 
(or by a loop). 
\end{remark}

\begin{example} 
Let $\cG$  consists of the following moves 
$(s, v)$, $(s, t)$, $(v, t')$, $(v, t'')$. 
Furthermore, $s$ is an initial position, 
while $t',t''$, and $t$ are terminal.  
The terminal payoffs are 
(2,3) in  $t$,  (3,2) in  $t'$, and  (1,1)  in $t''$. 
It is easy to verify that this game has two NE. 
The standard ``backward induction" NE:  
$x' = \{(s,v), (v, t')\}$ 
resulting in (3,2) and the blackmail  NE:  
$x'' = \{(s, t), (v, t'')$ resulting in (2,3). 
Player 2 threatens by move  $(v, t'')$, 
which is the worst for both, thus, pushing player 1 
to move to  $t$  which is better for  2  than  $t'$. 
Consider two terminal moves $(v, t')$  and $(v, t'')$  of player 2.  
The first is better than the second. 
Delete the bad move  $(v, t'')$  from the digraph. 
By doing so, we destroy the second NE. 
Only the first NE remains. 
\end{example}

\section{Main example}
The following example disproves CNJ$(C_{22})$.
However, 3 weaker conjectures, CNJ$(C'_{22})$,  
CNJ$(C_0)$ and CNJ$(C'_0)$ remain open. 

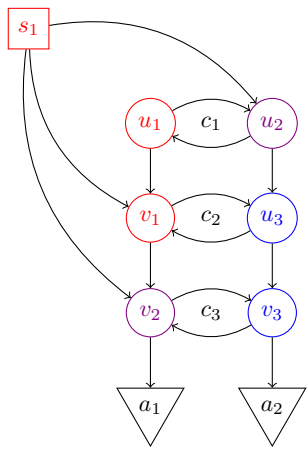
\begin{figure}[!ht]
\begin{center}
   \begin{tikzpicture}[scale=0.9, transform shape,node distance=14mm and 18mm, rnode/.style={draw,rectangle},mynode/.style={draw,circle}, pnode/.style={draw,regular polygon}, on grid]
        \node (s) [rnode, minimum size=6mm, red] {$s_1$};
        \node (u1) [mynode, below right=of s, red] {$u_1$};
        \node (u2) [mynode, right=of u1, violet] {$u_2$};
        \node (v1) [mynode, below=of u1, red] {$v_1$};
        \node (v2) [mynode, below=of u2, blue] {$u_3$};
        \node (w1) [mynode, below=of v1, violet] {$v_2$};
        \node (w2) [mynode, below = of v2, blue] {$v_3$};
        \node (a1) [pnode, regular polygon sides = 3, inner sep = 0.3mm, shape border rotate=180, below=of w1] {$a_1$};
        \node (a2) [pnode, regular polygon sides = 3, inner sep = 0.3mm, shape border rotate=180, below=of w2] {$a_2$};
        \graph [use existing nodes] {
            s -> [bend left] u2;
            u1 -> [bend left] u2; 
            u2 -> [bend left] u1;
            s -> [bend right] v1;
            s -> [bend right] w1;
            v1 -> [bend left] v2;
            v2 -> [bend left] v1;
            w1 -> [bend left] w2;
            w2 -> [bend left] w1;
            u1 -> v1 -> w1 -> a1;
            u2 -> v2 -> w2 -> a2;
        };
        \path(u1)--node{$c_1$}(u2);
        \node at ($(v1)!.5!(v2)$){$c_2$}; %requires calc library
        \node[right=9mm of w1] {$c_3$};
    \end{tikzpicture}
\end{center}
\caption{Graphical form}
\label{f1}
\end{figure}

The digraph $G = (V,E)$  is given in Figure \ref{f1}.
It has 3 distinct cyclic outcomes $c_1,c_2,c_3$ and only 
2 terminals $a_1, a_2$. 
The preferences of the players are also given Figure \ref{f1}. 
We see that $a_1$ and $a_2$ are top outcomes 
of players 3 and 1, respectively 
and, thus, the rank vector is $r = (1,2,1)$ 
and condition CND$(C_{22})$ fails. 
Nevertheless, we will show that 
the game is NE-free, in contradiction with CNJ$(C_{22})$. 

Three players $I = \{1,2,3\}$  control 
positions $(s_1,u_1,v_1)$, $(u_2,v_2)$, $(u_3,v_3)$, respectively, 
and have $|X_1| = 3 \times 2 \times 2 = 12$ and  
$|X_2| = |X_3| = 2 \times 2 = 4$ strategies. 
The initial position is $s_1$.
To specify a strategy of a player $i$ 
we fix how (s)he moves in each position. 
The normal form is given in Figure \ref{f2}. 
For every situation $x \in X$ 
we provide the outcome 
$\go(x) \in \gO = \{a_1,a_2,c_1,c_2,c_3\}$ 
with the superscripts showing each player $i$ 
who can improve $x$ by replacing  
$x_i$ with $x'_i$ such that  $\go(x') \succ_i g(x)$. 
In some situations there are more than one superscript, 
yet, for each situation there exists at least one.
This exactly means that the obtained game is NE-free. 
Yet, condition CND$(C_{22})$ of the conjecture CNJ$(C_{22})$ fail. 

However, if we merge cyclic outcomes 
$c_1,c_2,c_3$ then an NE will appear. 
Thus, a weaker conjecture CNJ$(C'_{22})$ 
is not disproved, as well as conjectures
CNJ$(C_0)$ and CNJ$(C'_0)$, since the corresponding 
conditions CND$(C_0)$ and CND$(C'_0)$ do not hold.

\renewcommand{\arraystretch}{1,5}

%The payoff functions $f_i$:
%$$ \text{player 1: } f_1(c_1) < f_1(a_1) < f_1(c_2) < f_1(c_3) < f_1(a_2) $$
%$$ \text{player 2: } f_2(a_2) < f_2(c_1) < f_2(c_2) < f_2(a_1) < f_2(c_3) $$
%$$ \text{player 3: } f_3(c_3) < f_3(c_1) < f_3(a_2) < f_3(c_2) < f_3(a_1) $$

Preferences $\prec_1, \prec_2, \prec_3$:
$$ \text{player 1: } c_1 \prec_1 a_1 \prec_1 c_2 \prec_1 c_3 \prec_1 a_2 $$
$$ \text{player 2: } a_2 \prec_2 c_1 \prec_2 c_2 \prec_2 a_1 \prec_2 c_3 $$
$$ \text{player 3: } c_3 \prec_3 c_1 \prec_3 a_2 \prec_3 c_2 \prec_3 a_1 $$

\clearpage

\begin{figure}[H]
        \centering
        \begin{subfigure}[!ht]{\textwidth}
           3-d player: $(v_2 \to v_1, w_2 \to w_1)$
%1. The third player's strategy is 

    \begin{center}
\scriptsize
\begin{tblr}{colspec = {*{5}{c}},
             cell{1}{1} ={r=2}{},
             cell{1}{2} ={c=4}{}
             }
  \hline
    \toprule
1st player    &   2nd player
                    &   &   &                               \\
    \midrule
                &   {$(u_2 \to u_1,$\\ $v_2 \to v_3)$}   
                    &   {$(u_2 \to u_1,$\\ $v_2 \to a_1)$}    
                        &   {$(u_2 \to u_3,$\\ $v_2 \to v_3)$}  
                            &   {$(u_2 \to u_3,$\\ $v_2 \to a_1)$}  \\  
    \midrule
\hline
$(s \to u_2, u_1 \to u_2, v_1 \to v_2)$ & $c_1^{1, 2}$ & $c_1^{1, 2}$ & $c_2^{1}$ & $c_2^{3}$ \\ 
$(s \to u_2, u_1 \to u_2, v_1 \to w_1)$ & $c_1^{1, 2}$ & $c_1^{1, 2}$ & $c_3^{3}$ & $a_2^{1, 2}$ \\ 
$(s \to u_2, u_1 \to v_1, v_1 \to v_2)$ & $c_2^{1}$ & $c_2^{3}$ & $c_2^{1}$ & $c_2^{3}$ \\ 
$(s \to u_2, u_1 \to v_1, v_1 \to w_1)$ & $c_3^{3}$ & $a_2^{1, 2}$ & $c_3^{3}$ & $a_2^{1, 2}$ \\ 
$(s \to v_1, u_1 \to u_2, v_1 \to v_2)$ & $c_2^{1}$ & $c_2^{3}$ & $c_2^{1}$ & $c_2^{3}$ \\ 
$(s \to v_1, u_1 \to u_2, v_1 \to w_1)$ & $c_3^{3}$ & $a_2^{1, 2}$ & $c_3^{3}$ & $a_2^{1, 2}$ \\ 
$(s \to v_1, u_1 \to v_1, v_1 \to v_2)$ & $c_2^{1}$ & $c_2^{3}$ & $c_2^{1}$ & $c_2^{3}$ \\ 
$(s \to v_1, u_1 \to v_1, v_1 \to w_1)$ & $c_3^{3}$ & $a_2^{1, 2}$ & $c_3^{3}$ & $a_2^{1, 2}$ \\ 
$(s \to w_1, u_1 \to u_2, v_1 \to v_2)$ & $c_3^{3}$ & $a_2^{1, 2}$ & $c_3^{3}$ & $a_2^{1, 2}$ \\ 
$(s \to w_1, u_1 \to u_2, v_1 \to w_1)$ & $c_3^{3}$ & $a_2^{1, 2}$ & $c_3^{3}$ & $a_2^{1, 2}$ \\ 
$(s \to w_1, u_1 \to v_1, v_1 \to v_2)$ & $c_3^{3}$ & $a_2^{1, 2}$ & $c_3^{3}$ & $a_2^{1, 2}$ \\ 
$(s \to w_1, u_1 \to v_1, v_1 \to w_1)$ & $c_3^{3}$ & $a_2^{1, 2}$ & $c_3^{3}$ & $a_2^{1, 2}$ \\ [1ex] 
    \bottomrule
\end{tblr}
    \end{center} 
        \end{subfigure}

        \begin{subfigure}[H]{\textwidth}
      \vspace{5mm}   
     3-d player: $(v_2 \to v_1, w_2 \to a_2)$
% 2. The third player's strategy is     
       \begin{center}
\scriptsize
\begin{tblr}{colspec = {*{5}{c}},
             cell{1}{1} ={r=2}{},
             cell{1}{2} ={c=4}{}
             }
  \hline
    \toprule
1st player    &   2nd player
                    &   &   &                               \\
    \midrule
                &   {$(u_2 \to u_1,$\\ $w_1 \to w_2)$}   
                    &  { $(u_2 \to u_1,$\\ $w_1 \to a_1)$}    
                        &   {$(u_2 \to v_2,$\\ $w_1 \to w_2)$}  
                            &   {$(u_2 \to v_2,$\\ $w_1 \to a_1)$}  \\  
    \midrule
\hline
$(s \to u_2, u_1 \to u_2, v_1 \to v_2)$ & $c_1^{1, 2}$ & $c_1^{1, 2}$ & $c_2^{1}$ & $c_2^{3}$ \\ 
$(s \to u_2, u_1 \to u_2, v_1 \to w_1)$ & $c_1^{1, 2}$ & $c_1^{1, 2}$ & $a_1^{2}$ & $a_2^{1}$ \\ 
$(s \to u_2, u_1 \to v_1, v_1 \to v_2)$ & $c_2^{1}$ & $c_2^{3}$ & $c_2^{1}$ & $c_2^{3}$ \\ 
$(s \to u_2, u_1 \to v_1, v_1 \to w_1)$ & $a_1^{2}$ & $a_2^{1}$ & $a_1^{2}$ & $a_2^{1}$ \\ 
$(s \to v_1, u_1 \to u_2, v_1 \to v_2)$ & $c_2^{1}$ & $c_2^{3}$ & $c_2^{1}$ & $c_2^{3}$ \\ 
$(s \to v_1, u_1 \to u_2, v_1 \to w_1)$ & $a_1^{2}$ & $a_2^{1}$ & $a_1^{2}$ & $a_2^{1}$ \\ 
$(s \to v_1, u_1 \to v_1, v_1 \to v_2)$ & $c_2^{1}$ & $c_2^{3}$ & $c_2^{1}$ & $c_2^{3}$ \\ 
$(s \to v_1, u_1 \to v_1, v_1 \to w_1)$ & $a_1^{2}$ & $a_2^{1}$ & $a_1^{2}$ & $a_2^{1}$ \\ 
$(s \to w_1, u_1 \to u_2, v_1 \to v_2)$ & $a_1^{2}$ & $a_2^{1}$ & $a_1^{2}$ & $a_2^{1}$ \\ 
$(s \to w_1, u_1 \to u_2, v_1 \to w_1)$ & $a_1^{2}$ & $a_2^{1}$ & $a_1^{2}$ & $a_2^{1}$ \\ 
$(s \to w_1, u_1 \to v_1, v_1 \to v_2)$ & $a_1^{2}$ & $a_2^{1}$ & $a_1^{2}$ & $a_2^{1}$ \\ 
$(s \to w_1, u_1 \to v_1, v_1 \to w_1)$ & $a_1^{2}$ & $a_2^{1}$ & $a_1^{2}$ & $a_2^{1}$ \\ [1ex] 
    \bottomrule
\end{tblr}
    \end{center}       
        \end{subfigure}

    \end{figure}

    \begin{figure}[H]\ContinuedFloat
        \centering

        \begin{subfigure}[h]{\textwidth}
    3-d player: $(v_2 \to w_2, w_2 \to w_1)$
%3. The third player's strategy is 

    \begin{center}
\scriptsize
\begin{tblr}{colspec = {*{5}{c}},
             cell{1}{1} ={r=2}{},
             cell{1}{2} ={c=4}{}
             }
  \hline
    \toprule
1st player    &   2nd player
                    &   &   &                               \\
    \midrule
                &   {$(u_2 \to u_1,$\\ $w_1 \to w_2)$}   
                    &  { $(u_2 \to u_1,$\\ $w_1 \to a_1)$}    
                        &   {$(u_2 \to v_2,$\\ $w_1 \to w_2)$}  
                            &   {$(u_2 \to v_2,$\\ $w_1 \to a_1)$}  \\  
    \midrule
\hline
$(s \to u_2, u_1 \to u_2, v_1 \to v_2)$ & $c_1^{1, 2}$ & $c_1^{1, 2}$ & $c_3^{3}$ & $a_2^{2}$ \\ 
$(s \to u_2, u_1 \to u_2, v_1 \to w_1)$ & $c_1^{1, 2}$ & $c_1^{1, 2}$ & $c_3^{3}$ & $a_2^{2}$ \\ 
$(s \to u_2, u_1 \to v_1, v_1 \to v_2)$ & $c_3^{3}$ & $a_2^{2}$ & $c_3^{3}$ & $a_2^{2}$ \\ 
$(s \to u_2, u_1 \to v_1, v_1 \to w_1)$ & $c_3^{3}$ & $a_2^{2}$ & $c_3^{3}$ & $a_2^{2}$ \\ 
$(s \to v_1, u_1 \to u_2, v_1 \to v_2)$ & $c_3^{3}$ & $a_2^{2}$ & $c_3^{3}$ & $a_2^{2}$ \\ 
$(s \to v_1, u_1 \to u_2, v_1 \to w_1)$ & $c_3^{3}$ & $a_2^{2}$ & $c_3^{3}$ & $a_2^{2}$ \\ 
$(s \to v_1, u_1 \to v_1, v_1 \to v_2)$ & $c_3^{3}$ & $a_2^{2}$ & $c_3^{3}$ & $a_2^{2}$ \\ 
$(s \to v_1, u_1 \to v_1, v_1 \to w_1)$ & $c_3^{3}$ & $a_2^{2}$ & $c_3^{3}$ & $a_2^{2}$ \\ 
$(s \to w_1, u_1 \to u_2, v_1 \to v_2)$ & $c_3^{3}$ & $a_2^{2}$ & $c_3^{3}$ & $a_2^{2}$ \\ 
$(s \to w_1, u_1 \to u_2, v_1 \to w_1)$ & $c_3^{3}$ & $a_2^{2}$ & $c_3^{3}$ & $a_2^{2}$ \\ 
$(s \to w_1, u_1 \to v_1, v_1 \to v_2)$ & $c_3^{3}$ & $a_2^{2}$ & $c_3^{3}$ & $a_2^{2}$ \\ 
$(s \to w_1, u_1 \to v_1, v_1 \to w_1)$ & $c_3^{3}$ & $a_2^{2}$ & $c_3^{3}$ & $a_2^{2}$ \\ [1ex] 
    \bottomrule
\end{tblr}
    \end{center}        
        \end{subfigure}
        \begin{subfigure}[H]{\textwidth}
     \vspace{5mm}   
     3-d player:  $(v_2 \to w_2, w_2 \to a_2)$ 
    \begin{center}
\scriptsize
\begin{tblr}{colspec = {*{5}{c}},
             cell{1}{1} ={r=2}{},
             cell{1}{2} ={c=4}{}
             }
  \hline
    \toprule
1st player    &   2nd player
                    &   &   &                               \\
    \midrule
                &   {$(u_2 \to u_1,$\\ $w_1 \to w_2)$}   
                    &  { $(u_2 \to u_1,$\\ $w_1 \to a_1)$}    
                        &   {$(u_2 \to v_2,$\\ $w_1 \to w_2)$}  
                            &   {$(u_2 \to v_2,$\\ $w_1 \to a_1)$}  \\  
    \midrule
\hline
$(s \to u_2, u_1 \to u_2, v_1 \to v_2)$ & $c_1^{1}$ & $c_1^{1}$ & $a_1^{2, 3}$ & $a_1^{2, 3}$ \\ 
$(s \to u_2, u_1 \to u_2, v_1 \to w_1)$ & $c_1^{1}$ & $c_1^{1}$ & $a_1^{2}$ & $a_1^{2, 3}$ \\ 
$(s \to u_2, u_1 \to v_1, v_1 \to v_2)$ & $a_1^{3}$ & $a_1^{3}$ & $a_1^{3}$ & $a_1^{3}$ \\ 
$(s \to u_2, u_1 \to v_1, v_1 \to w_1)$ & $a_1^{2}$ & $a_2^{1}$ & $a_1^{2}$ & $a_1^{2, 3}$ \\ 
$(s \to v_1, u_1 \to u_2, v_1 \to v_2)$ & $a_1^{3}$ & $a_1^{3}$ & $a_1^{3}$ & $a_1^{3}$ \\ 
$(s \to v_1, u_1 \to u_2, v_1 \to w_1)$ & $a_1^{2}$ & $a_2^{1}$ & $a_1^{2}$ & $a_2^{1}$ \\ 
$(s \to v_1, u_1 \to v_1, v_1 \to v_2)$ & $a_1^{3}$ & $a_1^{3}$ & $a_1^{3}$ & $a_1^{3}$ \\ 
$(s \to v_1, u_1 \to v_1, v_1 \to w_1)$ & $a_1^{2}$ & $a_2^{1}$ & $a_1^{2}$ & $a_2^{1}$ \\ 
$(s \to w_1, u_1 \to u_2, v_1 \to v_2)$ & $a_1^{2}$ & $a_2^{1}$ & $a_1^{2}$ & $a_2^{1}$ \\ 
$(s \to w_1, u_1 \to u_2, v_1 \to w_1)$ & $a_1^{2}$ & $a_2^{1}$ & $a_1^{2}$ & $a_2^{1}$ \\ 
$(s \to w_1, u_1 \to v_1, v_1 \to v_2)$ & $a_1^{2}$ & $a_2^{1}$ & $a_1^{2}$ & $a_2^{1}$ \\ 
$(s \to w_1, u_1 \to v_1, v_1 \to w_1)$ & $a_1^{2}$ & $a_2^{1}$ & $a_1^{2}$ & $a_2^{1}$ \\ [1ex] 
    \bottomrule
\end{tblr}
    \end{center}      
        \end{subfigure}
        \caption{Normal form}
        \label{f2}
    \end{figure}

\section{NE-free DG $n$-person games and satisfiability} 
Constructing NE-free $n$-person DG games requires a computer. 
The problem is reduced to satisfiability and solved by a SAT-solver. 
The CNF corresponding to a DG game is "bipolar": 
part of its clauses are positive, 
no negated variables, while the remaining clauses 
are negative, all variables are negated. 

Introduce $\frac{1}{2}nk(k-1)$ variables 
$y_{\go, \go'}^i$,   
for all $i \in I = \{1, \dots, n\}$ 
and distinct pairs $\go, \go' \in \gO$; 
here $k = |\gO|$. 
We set  $y_{\go, \go'}^i = 1$ when 
$\go \succ_i \go'$ and, respectively, 
$y_{\go, \go'}^i = 0$ when  $\go' \succ_i \go$. 
Thus, $y_{\go, \go'}^i = \bar{y}_{\go', \go}^i$ 
holds by definition. Consider CNF 
$C_- = \bigwedge (\bar{y}_{\go, \go'}^i \vee 
 \bar{y}_{\go', \go''}^i \vee 
 \bar{y}_{\go'', \go}^i)$, 
 where conjunction is taken over  
 all players $i \in I$  and all 
 pairwise distinct triplets 
 $\go, \go', \go'' \in \gO$. 
 % By definition, 
 All variables of  $C_-$ are negated.
 It is not difficult to see that $C_-$  is satisfied 
 if and only if it corresponds to an acyclic 
 preference of player $i$  over $\gO$. 
 Indeed, 2-dicycles obviously are not possible, 
 since  $y_{\go, \go'}^i = \bar{y}_{\go', \go}^i$. 
 Furthermore, it is well-known that 
 any tournament that has a dicycle,  
 has a 3-dicycle; see, for example, \cite{Tho86}. 
Finally,  it is not difficult to 
assign a positive CNF $C_+$ to a game such that 
the game is NE-free if and only if 
CNF $C_i \wedge C_+$ is satisfyable. 
Although it is exponential in the size of digraph  $G$, 
but its structure is clear. 
Then, by the above construction, 
an NE-free preference $\succ$ 
for the considered DG $n$-person game 
structure corresponds to a satisfying assignments 
of the CNF $C_- \wedge C_+$. 
Using a powerful SAT-solver we can analyze  
many games and find NE-free ones among them.

\section{NE-free $n$-person DG games 
with a unique terminal outcome and 
the rank vector $r = (0,\dots,0,1)$} 

Consider an arbitrary NE-free $n$-person game $\Gamma$, 
% (with $n \geq 3)$. In particular, take the 
for example, one from \cite{BGMOV18} or from the present paper. 
The following simple trick allows us to construct 
a NE-free $n$-person game $\Gamma$ 
with the rank vector $r = (0,\dots,0,1)$. 
To do so, to each terminals $v \in V_T$ of $\Gamma$ 
add a loop  $l_v$  and assign the same player, say  $n$, 
and the same preferences as for the terminal $v$ in $\Gamma$.
Then add a new terminal $v^*$ to $\Gamma$ and 
a move from each $v \in V_T$  to $v^*$. 
Furthermore, let the unique new terminal
$v^*$ will be the worst outcome for player $n$ 
and the best outcome for each of the remaining $n-1$ players.
Obviously, the obtained game $\Gamma'$ is NE-free, 
it still has $n$ players, a unique terminal outcome $v^*$,  
and the rank vector $r = (0,\dots,0,1)$. 

Note finally that $\Gamma'$ has several non-terminal SCCs 
and also that in case of a unique 
non-terminal outcome 
(that is, when all non-terminal SCCs are merged) 
Conjecture Catch 22 remains open. 

\section{On NE-free play-once 
$n$-person DG games 
with a unique terminal outcome and a (0,1) rank vector} 

Currently, no NE-free play-once $n$-person DG game is known. 
However, any such example with $m$ terminal outcomes 
could be easily transformed to 
a NE-free play-once $(n+m)$-person DG game 
with a unique terminal outcome and 
a (0,1) rank vector that has   
$m$ entries 1 and $n$ entries  0. 

The transformation of the previous section works  
with the following minor modification. 
To each of the $m$ terminal positions of the original game 
we assign a separate new player 
(while  in the previous section 
a single player $n$ was assigned to all $m$ terminals). 
Furthermore, let each them prefer  
his/her loop to the unique new terminal outcome. 
Otherwise, the latter is better than 
any other outcome for each player. 
Obviously, this construction results in 
a play-once NE-free $(n+m)$ person DG game 
whose rank vector has $m$ entries 1 and $n$ entries 0.

\subsection*{Acknowledgements}
The paper was prepared within the framework
of the HSE University Basic Research Program.
% Moscow, Russia.
% and funded by the RSF grant  20-11-20203.
% The author is thankful to Endre Boros for many helpful remarks.

\end{document}